\documentclass[12pt,a4paper]{article}
\usepackage{amsmath,amsfonts,amssymb,amsthm}

\allowdisplaybreaks[1]

\newcommand{\J}{\mu_1}
\newcommand{\M}{\mu_2}
\newcommand{\trace}{\operatorname{trace}}

\newcommand{\Span}{\mathrm{span}}
\newcommand{\R}{\mathbb{R}}
\newcommand{\Rf}{\R^4}

\newcommand{\ba}{\begin{array}}
\newcommand{\ea}{\end{array}}

\newcommand{\iz}{\varphi}
\newcommand{\ic}[3]{{\iz_{#1 #2}^{#3}}}

\newcommand{\pt}{\frac{\partial}{\partial t}}

\newcommand{\half}{{\tfrac{1}{2}}}

\newcommand{\ricm}{\widehat{\operatorname{Ric}}}

\newtheorem{thm}{Theorem}
\newtheorem{lem}{Lemma}
\newtheorem*{question}{Question}

\theoremstyle{remark}
\newtheorem*{rem}{Remark}

\newcounter{rom}
\renewcommand{\therom}{(\roman{rom})}
{\end{list}}
\begin{document}

\title{3-dimensional affine hypersurfaces admitting a pointwise
  $SO(2)$- or $\mathbb Z_3$-symmetry} 
\author{C. Scharlach\thanks{Partially
    supported by the DFG-project SI 163-7} \and L.
  Vrancken\thanks{Partially supported by a research fellowship of the
    Alexander von Humboldt Stiftung (Germany)}} 
\maketitle

\begin{abstract} In (equi-)affine differential geometry, the most important
  algebraic invariants are the affine (Blaschke) metric $h$, the
  affine shape operator $S$ and the difference tensor $K$. A
  hypersurface is said to admit a pointwise symmetry if at every point
  there exists a linear transformation preserving the affine metric,
  the affine shape operator and the difference tensor $K$.  In this
  paper, we consider the $3$-dimensional positive definite
  hypersurfaces for which at each point the group of symmetries is
  isomorphic to either $\mathbb Z_3$ or $SO(2)$. We classify such
  hypersurfaces and show how they can be constructed starting from
  $2$-dimensional positive definite affine spheres.
 \end{abstract}

\medskip\noindent
{\bfseries Subject class: } 53A15


\section{Introduction}\label{sec:intro}
In this paper we study nondegenerate (equi-)affine hypersurfaces
$F\colon M^n \to \R^{n+1}$. In that case, it is well known that there exists
a canonical choice of transversal vector field $\xi$ called the affine
(Blaschke) normal, which induces a connection $\nabla$, a symmetric bilinear
form $h$ and a 1-1 tensor field $S$ by
\begin{align} 
&D_X Y =\nabla_X Y +h(X,Y)\xi,\label{strGauss}\\
&D_X \xi =-SX,\label{strWeingarten}
\end{align}
for all $X,Y \in {\cal X}(M)$. The connection $\nabla$ is called the induced
affine connection, $h$ is called the affine metric (or Blaschke
metric) and $S$ is called the affine shape operator.  In general $\nabla$
is not the Levi Civita connection $\hat\nabla$ of $h$. The difference
tensor $K$ is defined as
\begin{equation}\label{defK}
K(X,Y)=\nabla_X Y-\hat\nabla_X Y,
\end{equation}
for all $X,Y \in {\cal X}(M)$. Moreover the form $h(K(X,Y),Z)$ is a
symmetric cubic form with the property that for any fixed $X\in {\cal
  X}(M)$, $\trace K_X$ vanishes.  This last property is called the
apolarity condition.  The difference tensor $K$, together with the
affine metric $h$ and the affine shape operator are the most
fundamental algebraic invariants for a nondegenerate affine
hypersurface (more details in Sec.~\ref{sec:basics}). We say that $M$
is positive definite if the affine metric $h$ is positive definite.
For the basic theory of nondegenerate affine hypersurfaces we refer to
\cite{LSZ} and \cite{NS}.

A hypersurface is said to admit a pointwise symmetry if at every point
there exists a linear transformation preserving the affine metric, the
affine shape operator and the difference tensor $K$.  The study of
submanifolds which admit pointwise isometries was initiated by Bryant
in \cite{B} where he studied $3$-dimensional Lagrangian submanifolds
of $\mathbb C^3$. Following essentially the same approach, a
classification of $3$-dimensional affine hyperspheres admitting
pointwise isometries was obtained in \cite{V}.

In \cite{sv}, for $3$-dimensional positive definite hypersurfaces, the
possible groups which can act on the algebraic invariants as well as
the canonical forms for $S$, $K$ and $h$ were computed.  In this
paper, we consider the $3$-dimensional positive definite hypersurfaces
for which at each point the group of symmetries is isomorphic to
either $\mathbb Z_3$ or the group of rotations $SO(2)$.  The paper is
organized as follows. First in Section 2, we shortly recall the basic
equations of Gauss, Codazzi and Ricci for an affine hypersurface and
use those equations, together with the canonical form of $h$, $S$ and
$K$, to obtain information about the coefficients of the connection.
In particular, it follows that such a hypersurface $M$ admits a warped
product structure.  In Section 3, we classify such hypersurfaces by
showing how they can be constructed starting from $2$-dimensional
positive definite affine spheres. This classification can be seen as a
generalisation of the well known Calabi product of hyperbolic affine
spheres and of the constructions for affine spheres considered in
\cite{DV}. Note that affine hyperspheres, i.e. affine hypersurfaces
for which all affine normals are parallel or pass through a fixed
point, are without any doubt the most studied class of affine
hypersurfaces.  They are closely related to solutions of Monge Amp{\`e}re
equations.  The following natural question for a (de)composition
theorem, related to the Calabi product and its generalisations in
\cite{DV}, gives another motivation for studying $3$-dimensional
hypersurfaces admitting a $\mathbb Z_3$-symmetry or an $SO(2)$-symmetry:

\begin{question} Let $M^n$ be a nondegenerate affine hypersurface in
  $\mathbb R^{n+1}$. Under what conditions do there exist affine
  hyperpsheres $M_1^r$ in $\mathbb R^{r+1}$ and $M_2^s$ in $\mathbb
  R^{s+1}$, with $r+s=n-1$, such that $M = I \times_{f_1} M_1 \times_{f_2} M_2$, where
  $I \subset \mathbb R$ and $f_1$ and $f_2$ depend only on $I$ (i.e. $M$
  admits a warped product structure)? How can the original immersion
  be recovered starting from the immersion of the affine spheres?
\end{question}

Of course the first dimension in which the above problem can be
considered is three and our study of $3$-dimensional affine
hypersurfaces with $\mathbb Z_3$-symmetry or $SO(2)$-symmetry
provides an answer in that case.

\section{Structure equations and integrability conditions}

\subsection{Preliminaries}
\label{sec:basics}

We consider $3$-dimensional affine hypersurfaces $F\colon M^3 \to\Rf$.
Assume that $M^3$ has at every point a $\mathbb Z_3$-symmetry or an
$SO(2)$-symmetry.  Then we recall from \cite{sv} the following: At every
point $p$ of $M^3$ there exists a basis $\{e_1,e_2,e_3\}$ which is
orthonormal with respect to the affine metric $h$ such that the
difference tensor $K$ and the shape operator $S$ are respectively
given by:
 \begin{equation}\label{pdifftensor}
\begin{split}
  K_{e_1}&=\begin{pmatrix} 2\J &0&0\\ 0& -\J & 0 \\ 0&0& -\J
  \end{pmatrix}, \quad
  K_{e_2}=\begin{pmatrix} 0 & -\J&0\\ -\J & \M & 0 \\ 0&0& -\M
  \end{pmatrix}, \\
K_{e_3}&=\begin{pmatrix} 0 & 0 & -\J\\ 0 & 0 & -\M \\ -\J &-\M& 0
  \end{pmatrix},\quad 
S=\begin{pmatrix} \lambda& 0 &0 \\ 0& a &0 \\ 0& 0 & a
  \end{pmatrix}.
\end{split} \end{equation}
We have that $\J$ is nonzero. Moreover, $\M$ vanishes if and only
if $M^3$ admits a $1$-parameter group $SO(2)$ of isometries. In that case
the form of $K$ and $S$ remains invariant under rotations in the
$e_2e_3$-plane.  In case that $\M$ is different from zero, the group
$\mathbb Z_3$ of rotations leaving $K$ and $S$ invariant is generated
by the rotation with angle $\tfrac{2\pi}{3}$ in the $e_2e_3$-plane.

We recall some of the fundamental equations, which a nondegenerate
hypersurface has to satisfy, see also \cite{NS} or \cite{LSZ}. These
equations relate $S$ and $K$ with amongst others the curvature tensor
$R$ of the induced connection $\nabla$ and the curvature tensor $\hat
R$ of the Levi Civita connection $\hat\nabla$ of the affine metric
$h$.  We respectively have the Gauss equation for $\nabla$, which
states that:
\begin{equation}\label{gauss}
R(X,Y)Z =h(Y,Z)SX -h(X,Z) SY,
\end{equation}
and the Codazzi equation 
\begin{equation}\label{codazzi}
(\nabla_X S) Y =(\nabla_Y S) X.
\end{equation}
The fundamental existence and uniqueness theorem, see \cite{D} or
\cite{DNV}, states that given $h$, $\nabla$ and $S$ such that the
difference tensor is symmetric and traceless with respect to $h$, on a
simply connected manifold $M$ an affine immersion of $M$ exists if and
only if the above Gauss equation and Codazzi equation are satisfied.
From these the Codazzi equation for $K$ and the Gauss equation for
$\hat\nabla$ follow.
\begin{equation}\label{CodK2}\begin{split}
(\hat{\nabla}_X K)(Y,Z)- (\hat{\nabla}_Y K)(X,Z)=& \half( h(Y,Z)SX-
h(X,Z)SY\\& - h(SY,Z)X + h(SX, Z)Y),\end{split}
\end{equation}
and
\begin{equation}\label{gaussmetric}\begin{split} \hat R (X,Y)Z =&
\tfrac{1}{2} (h(Y,Z) SX - h(X,Z) SY\\& + h(SY,Z)X - h(SX,Z)Y ) -
    [K_X ,K_Y] Z\end{split}\end{equation}

\subsection{An adapted frame}
\label{sec:frame}

From now on, we assume that $M^3$ admits a $\mathbb Z_3$-symmetry or
an $SO(2)$-symmetry. The first meaning that at every point of $M^3$ the
group of isometries preserving $S$ and $K$ is isomorphic to $\mathbb
Z_3$, whereas in the second case, we assume that that group of
isometries is at every point isomorphic to $SO(2)$.

We define the Ricci tensor of the connection $\hat \nabla$ by:
\begin{equation*}
\ricm(X,Y)=\trace\{ Z \mapsto \hat R(Z,X)Y\}.
\end{equation*}
It is well known that $\ricm$ is a symmetric operator.
Then, we have
\begin{lem} Let $p \in M$ and $\{e_1,e_2,e_3\}$ the basis constructed
  earlier. Then 
\begin{alignat*}{2}
  &\ricm(e_1,e_1) = (a+\lambda)+ 6 \J^2, \qquad\quad &&\ricm(e_1,e_2)=0, \\
  &\ricm(e_3,e_1)=0,\qquad\quad &&
\ricm(e_2,e_2)=\tfrac{3}{2}a +\tfrac{1}{2} \lambda +2(\J^2+\M^2),\\
  &\ricm(e_2,e_3)=0,&&\ricm(e_3,e_3)=\tfrac{3}{2}a +\tfrac{1}{2}
  \lambda +2(\J^2+\M^2).
\end{alignat*}
\end{lem}
\begin{proof}We use the Gauss equation \eqref{gaussmetric} for $\hat R$.
It follows that
\begin{align*}
  \hat R(e_2,e_1)e_1&=\tfrac{1}{2}(a+\lambda) e_2-K_{e_2}(2\J
  e_1)+K_{e_1}(-\J e_2) \\
  &=(\tfrac{1}{2}(a+\lambda)+3 \J^2) e_2,\\
  \hat R(e_3,e_1)e_1&=\tfrac{1}{2}(a+\lambda) e_3-K_{e_3}(2\J e_1)
  +K_{e_1}(-\J e_3) \\
  &=(\tfrac{1}{2}(a+\lambda)+3 \J^2) e_3,\\
  \hat R(e_3,e_1)e_2&=-K_{e_3}(-\J e_2)+K_{e_1}(-\M e_3)=0.
\end{align*}
From this it immediately follows that 
$$\ricm(e_1,e_1) = (a+\lambda)+6 \J^2$$
and
$$\ricm(e_1,e_2)=0.$$
The other equations follow by similar computations.
\end{proof}
Now, we want to show that the basis we have constructed at each point
$p$ can be extended differentiably to a neighborhood of the point $p$
such that at every point the components of $S$ and $K$ with respect to
the frame $\{e_1,e_2,e_3\}$ have the previously described form.
\begin{lem} Let $M^3$ be an affine hypersurface of $\mathbb R^4$ which
  admits a pointwise $\mathbb Z_3$-symmetry or a pointwise
  $SO(2)$-symmetry. Let $p \in M$. Then there exists a frame
  $\{e_1,e_2,e_3\}$ defined in a neighborhood of the point $p$ such that
  the components of $K$ and $S$ are respectively given by:
\begin{equation*}
\begin{split}
  K_{e_1}&=\begin{pmatrix} 2\J &0&0\\ 0& -\J & 0 \\ 0&0& -\J
  \end{pmatrix}, \quad
  K_{e_2}=\begin{pmatrix} 0 & -\J&0\\ -\J & \M & 0 \\ 0&0& -\M
  \end{pmatrix}, \\
K_{e_3}&=\begin{pmatrix} 0 & 0 & -\J\\ 0 & 0 & -\M \\ -\J &-\M& 0
  \end{pmatrix},\quad 
S=\begin{pmatrix} \lambda& 0 &0 \\ 0& a &0 \\ 0& 0 & a
  \end{pmatrix}.
\end{split} 
\end{equation*}
 \end{lem}
\begin{proof}First we want to show that at every point the vector
  $e_1$ is uniquely defined and differentiable. We introduce a
  symmetric operator $\hat A$ by:
\begin{equation*}
\ricm(Y,Z)= h(\hat A Y,Z).
\end{equation*}
Clearly $\hat A$ is a differentiable operator on $M$. On the set of
points where $-\tfrac{1}{2}(a-\lambda) +4\J^2 -2 \M^2\neq 0$, the
operator has two distinct eigenvalues. The eigendirection which
corresponds with the $1$-dimensional eigenvalue corresponds with the
vector field $e_1$. 

On the set of points where $-\tfrac{1}{2}(a-\lambda) +4\J^2 -2 \M^2=
0$, we have that $a \neq \lambda$ (cp. \cite{sv}), as otherwise we
would have an $A_4$-symmetry. In this case, the differentiable
operator $S$ has two distinct eigenvalues (cp. \eqref{pdifftensor})
and $e_1$ is uniquely determined as the eigendirection corresponding
to the $1$-dimensional eigenvalue. This shows that taking at every
point $p$ the vector $e_1$ yields a differentiable vector field. 

To show that $e_2$ and $e_3$ can be extended differentiably, we
consider two cases. First we assume that $M$ admits a pointwise
$SO(2)$-symmetry. In that case we have that $\M=0$ and we take for $e_2$
and $e_3$ arbitrary orthonormal differentiable local vector fields
which are orthogonal to the vector field $e_1$. In case that $M$
admits a pointwise $\mathbb Z_3$-symmetry we proceed as follows. We
start by taking arbitrary orthonormal differentiable local vector
fields $u_2$ and $u_3$ which are orthogonal to the vector field $e_1$.
It is then straightforward to check that we can write
\begin{align*}
  &K_{u_2}u_2= -\J e_1 +\nu_1 u_2 +\nu_2 u_3\\
  &K_{u_2}u_3= \nu_2 u_2 -\nu_1 u_3\\
  &K_{u_3}u_3= -\J e_1 -\nu_1 u_2 -\nu_2 u_3\end{align*} for some
differentiable functions $\nu_1$ and $\nu_2$ with $\nu_1^2
+\nu_2^2\neq 0$. Therefore, if necessary by interchanging the role of
$u_2$ and $u_3$, we may assume that in a neighborhood of the point
$p$, $\nu_1\neq 0$.  Rotating now over an angle $\theta$, thus
defining
\begin{align*}
&e_2 =\cos \theta u_2 +\sin \theta u_3,\\
&e_3 =-\sin \theta u_2 + \cos \theta u_3,
\end{align*}
we get that 
\begin{align*}
  h(K(e_2,e_2),e_2)&=(\cos^3 \theta -3 \cos \theta \sin^2 \theta)
  \nu_1 + (-\sin^3 \theta +3 \cos^2 \theta \sin \theta)\nu_2\\
  &=\cos 3 \theta \nu_1 +\sin 3 \theta \nu_2\\
  h(K(e_3,e_3),e_3)&=(-\sin^3 \theta +3 \cos^2 \theta \sin \theta)
  \nu_1 + (-\cos^3 \theta +3 \cos \theta \sin^2 \theta)\nu_2\\
  &=\sin 3 \theta \nu_1 - \cos 3 \theta \nu_2.
\end{align*}
Therefore, taking into account the symmetries of $K$, in order to
obtain the desired frame, it is sufficient to choose $\theta$ in such
a way that
\begin{equation*}
\sin 3 \theta \nu_1 - \cos 3 \theta \nu_2=0,
\end{equation*}
and $\cos 3 \theta \nu_1 +\sin 3 \theta \nu_2>0$. As this is always
possible, the proof is completed.
\end{proof}

\begin{rem} It actually follows from the proof of the previous lemma
  that the vector field $e_1$ is globally defined on $M$, and
  therefore the function $\J$, too. This in turn implies that the
  functions $\M$ (as it can be expressed in terms of $\J$ and the Pick
  invariant $J$ and $J$ is either identically zero or nowhere zero),
  $\lambda$ and $a$ (as it can be expressed in terms of the mean
  curvature and $\lambda$) are globally defined functions on the
  affine hypersurface $M$.
\end{rem}
From now on we will always work with the local frame constructed in
the previous lemma. We introduce the connection coefficients with
respect to this frame by $\hat \nabla_{e_i} e_j =\sum_{k=1}^3
\ic{i}{j}{k} e_k$. As the connection $\hat \nabla$ is metrical, we
have the usual symmetries.  

\subsection{Codazzi equations for $K$} \label{sec:CodK}

An evaluation of the Codazzi equations \eqref{CodK2} for $K$ using the
computer program mathematica (see also:
www-sfb288.math.tu-berlin.de/{\~{}}cs/symm2.nb resp. symm6.nb) results
in the following equations:
\begin{eqnarray}\allowdisplaybreaks
e_2(\J)=2 \J \ic112, & (eq. 1)\label{eq1.1}\\
\M \ic113 =4 \J \ic213, & (eq. 1)\label{eq1.2}\\
e_1(\J)= \tfrac{1}{2} (a-\lambda)-\M \ic112- 4\J \ic212,
&(eq. 1)\label{eq1.3}\\ 
e_1(\J)=\tfrac{1}{2} (a-\lambda)+\M \ic112 -4\J \ic313, & (eq. 2)
\label{eq2.1}\\
e_3(\J)=2 \J \ic113, & (eq. 2)\label{eq2.2}\\
\M \ic113=4 \J \ic312, & (eq. 2)\label{eq2.3}\\
e_1(\M)+e_2(\J)= 3\J \ic112 -\M \ic212, & (eq. 3)\label{eq3.1}\\
0=-\J \ic113 + 3\M \ic123 -\M \ic213, & (eq. 3)\label{eq3.2}\\
e_3(\J)= -\M (\ic213 +\ic312), & (eq. 4)\label{eq4.1}\\
e_3(\M)=3\M \ic223 -\J (\ic213-3\ic312), & (eq. 4)\label{eq4.2}\\
e_2(\M)=-\J (\ic212-\ic313) -3 \M \ic323, & (eq. 4)\label{eq4.3}\\ 
e_1(\M)=-\J\ic112 -\M \ic313, &(eq. 5)\label{eq5.1}\\
e_3(\J)=3 \J \ic113 +\M(3\ic123+\ic312), &(eq. 5)\label{eq5.2}\\
e_2(\J)=\M (\ic212-\ic313), &(eq. 6)\label{eq6.1}\\
e_3(\M)=3\M \ic223 +\J (3 \ic213-\ic312), & (eq. 6)\label{eq6.2}\\
e_2(\J)-e_1(\M)=\J \ic112 +\M \ic212, &(eq. 7)\label{eq7.1}\\
4 \J (\ic213 -\ic312) =0,&(eq. 8)\label{eq8.1}\\
e_3(\J)=\J \ic113- \M (3\ic123+\ic312). &(eq. 9) \label{eq9.1}
\end{eqnarray}

In the above expressions, the equation numbers refer to corresponding
equations in the mathematica program.  In order to simplify the above
equations, we now distinct two cases. 

\begin{lem}\label{evalCodK} An evaluatin of the Codazzi equations for
  $K$ gives:
\begin{gather*} \ic112=0,\quad \ic113=0,\quad \ic213=0,\quad \ic312=0,
  \quad \ic212= \ic313 =:\eta,\\
  e_1(\J)=\half(a-\lambda) -4 \J \eta,\quad e_2(\J)=0=e_3(\J).
\end{gather*}
If $\M\neq 0$, we get in addition that $\ic123=0$ and
\begin{equation}\label{diffmu2}
e_1(\M)= -\M \eta, \quad e_2(\M)= -3 \M \ic323,\quad 
e_3(\M)=3\M \ic223.
\end{equation}
\end{lem}

\begin{proof} First, we assume that $\M=0$ (thus $\J\neq 0$). In that
  case, it follows from \eqref{eq1.2} (resp. \eqref{eq2.3}) that
  $\ic213=0$ (resp.  $\ic312=0$), whereas \eqref{eq4.3} implies that
  $\ic212 =\ic313$. As it now follows from \eqref{eq4.1}and
  \eqref{eq6.1} that $e_2(\J)=e_3(\J)=0$, \eqref{eq1.1} and
  \eqref{eq2.2} imply that $\ic112=\ic113=0$. Finally \eqref{eq2.1}
  now reduces to
$$e_1(\J)= \tfrac{1}{2} (a-\lambda)- 4\J \ic212.$$

Next, we want to deal with the case that $\M \neq 0$. First it follows
from \eqref{eq1.1}, \eqref{eq6.1} and \eqref{eq3.1}, taking also into
account \eqref{eq5.1}, that
$$e_2(\J)=2 \J \ic112=\M (\ic212-\ic313)=4 \J \ic112 -\M
(\ic212-\ic313).$$
Therefore we get by \eqref{eq1.3} and
\eqref{eq2.1} that $\ic212=\ic313$ and thus $e_2(\J)=0=\ic112$.
From \eqref{eq4.1} and \eqref{eq2.2} it follows that
$$e_3(\J)=-\M (\ic213 +\ic312)=2 \J \ic113.$$
From \eqref{eq8.1},
\eqref{eq1.2} and the previous equation it follows that
$\ic213=\ic312=\ic113=0$ and $e_3(\J)=0$. From
\eqref{eq3.2} it then follows that $\ic123=0$.  From \eqref{eq3.1},
\eqref{eq4.3} and \eqref{eq4.2} we obtain the equations for $e_i(\M)$,
$i=1,2,3$ and from \eqref{eq1.3} it follows that
$$e_1(\J)= \tfrac{1}{2} (a-\lambda)- 4\J \ic212.$$
\end{proof}
As a direct consequence we write down the Levi-Civita connection:
\begin{lem}\label{LeviCivita}
\begin{eqnarray*}
   \hat{\nabla}_{e_1}e_1 &=& 0,\\
   \hat{\nabla}_{e_1}e_2 &=& \ic123 e_3,\\
   \hat{\nabla}_{e_1}e_3 &=& -\ic123 e_2,\\
   \hat{\nabla}_{e_2}e_1 &=& \eta e_2,\\
   \hat{\nabla}_{e_2}e_2 &=& -\eta e_1 + \ic223 e_3,\\
   \hat{\nabla}_{e_2}e_3 &=& -\ic223 e_2,\\
   \hat{\nabla}_{e_3}e_1 &=& \eta e_3,\\
   \hat{\nabla}_{e_3}e_2 &=& \ic323 e_3,\\
   \hat{\nabla}_{e_3}e_3 &=& -\eta e_1- \ic323 e_2,
 \end{eqnarray*}
where in case that $\M \neq 0$, we have in addition that
$\ic123=0$.
\end{lem}

\subsection{Gauss for $\nabla$} \label{sec:GaussIC}

Taking into account the previous results, we then proceed with an
evaluation of the Gauss equations \eqref{gauss} for $\nabla$:
\begin{equation*}
\nabla_X \nabla_Y Z- \nabla_Y \nabla_X Z - \nabla_{[X,Y]} Z = h(Y,Z)SX -
h(X,Z)SY,
\end{equation*}
again using the computer program mathematica (see also:
www-sfb288.math.tu-berlin.de/{\~{}}cs/symm2.nb resp. symm6.nb).  This
results amongst others in the following equations (cp. equations 11,
13, 14 and 16 in the mathematica program):
\begin{eqnarray*}
e_1(\eta)& =& -\eta^2 -3 \J^2 -\half(a+\lambda), \\
e_2(\eta)& =& 0,\\
e_3(\eta)& =&0.
\end{eqnarray*}

\subsection{Codazzi for $S$} \label{sec:CodS}

An evaluation of the Codazzi equations \eqref{codazzi} for $S$:
\begin{equation*}
(\nabla_X S)(Y) = (\nabla_Y S)(X)
\end{equation*}
by mathematica (see also:
www-sfb288.math.tu-berlin.de/{\~{}}cs/symm2.nb resp. symm6.nb,
equations 20 - 22)) then yields:
\begin{eqnarray*}
  e_1(a) &=& (\J-\eta)(a-\lambda),\\
e_2(a) &=& 0, \\
e_3(a) &=& 0, \\
e_2(\lambda) &=& 0, \\
e_3(\lambda) &=& 0.
\end{eqnarray*}

\subsection{Structure equations}
\label{sec:Streq} 

Summarized we have obtained the structure equations
(cp. \eqref{strGauss}, \eqref{strWeingarten} and \eqref{defK}):
\begin{alignat}{4}
&D_{e_1} e_1 =& 2\J e_1&&&+ \xi, \label{D11}\\
&D_{e_1} e_2 =& &-\J e_2 &+ \ic123 e_3,& \label{D12}\\
&D_{e_1} e_3 =& &-\ic123 e_2 &-\J e_3,& \label{D13}\\
&D_{e_2} e_1 =& &(\eta -\J) e_2, && \label{D21}\\
&D_{e_3} e_1 =& &&(\eta -\J) e_3, & \label{D31}\\
&D_{e_2} e_2 =&-(\eta+\J) e_1 &+\M e_2 &+ \ic223 e_3 &+\xi, \label{D22}\\
&D_{e_2} e_3 =& &-\ic223 e_2 &-\M e_3,& \label{D23}\\
&D_{e_3} e_2 =& &&(\ic323-\M) e_3, & \label{D32}\\
&D_{e_3} e_3 =&-(\eta+\J) e_1 &-(\ic323-\M) e_2 & &+\xi, \label{D33}
\end{alignat}\begin{alignat}{4}
&D_{e_1} \xi=& -\lambda e_1,&&&\label{D1xi}\\
&D_{e_2} \xi=&& -a e_2,&&\label{D2xi}\\
&D_{e_3} \xi=&&& -a e_3,&\label{D3xi}
\end{alignat}


Moreover, the functions $a$, $\lambda$, $\J$ and $\eta$ are all
constant in the $e_2$ and $e_3$-directions and the $e_1$-derivatives
are determined by (cp. Section~\ref{sec:CodS} and \ref{sec:GaussIC} and
Lemma~\ref{evalCodK}):
\begin{align}
&e_1(a) = (\J-\eta)(a-\lambda),\label{e1a}\\
&e_1(\eta) = -\eta^2 -3 \J^2 -\half(a+\lambda),\label{e1eta} \\
&e_1(\J)= - 4\J \eta-\tfrac{1}{2} (\lambda-a).\label{e1J}
\end{align}

\section{Main results}

As the vector field $e_1$ is globally defined, we can define the
distributions $H_1=\Span\{e_1\}$ and $H_2=\Span\{e_2,e_3\}$. In the
next lemmas we will investigate some properties of these distributions
following from Lemma~\ref{LeviCivita}. For the terminology we refer to
\cite{no}.
\begin{lem}\label{H1}
The distribution $H_1$ is autoparallel with respect to $\widehat\nabla$.
\end{lem}
\begin{proof} From $\hat{\nabla}_{e_1} e_1 =0$ the claim follows
  immediately. 
\end{proof} 

\begin{lem}\label{H2}
  The distribution $H_2$ is spherical with mean curvature normal
  $H=-\eta e_1$.
\end{lem}
\begin{proof} For $H=-\eta e_1\in H_1=H_2^{\perp}$ we have
  $h(\hat{\nabla}_{e_a} e_b, e_1)= h(e_a, e_b) h(H,e_1)$ for all $a,b
  \in \{2,3\}$, and $h(\hat{\nabla}_{e_a} H, e_1)= h(-e_a(\eta) e_1
  - \eta \hat{\nabla}_{e_a} e_1, e_1)=0$.
\end{proof}

\begin{rem} $\eta\; (=\ic212=\ic313)$ is independent of the
  choice of $e_2$ and $e_3$. It therefore is a globally defined
  function on $M$. 
\end{rem}

We introduce a coordinate function $t$ by $\pt:=e_1$. Using the
previous lemma, according to \cite{H}, we get:
\begin{lem}\label{warped} $(M,h)$ admits a warped product structure
  $M^3=\mathbb R\times_{e^f}N^2$ with $f:\mathbb R \to \mathbb R$
  satisfying
\begin{equation}\label{deff}
\frac{\partial f}{\partial t}=\eta.
\end{equation}
\end{lem}
\begin{rem} $a$, $\eta$ and $\J$ are functions of $t$, they satisfy
  by \eqref{e1a}, \eqref{e1eta} and \eqref{e1J}:
 \begin{align*}
   \frac{\partial a}{\partial t}&= (\J-\eta)(a-\lambda),\\
   \frac{\partial \eta}{\partial t}&=-\eta^2-3\J^2-\tfrac 12 (a+\lambda),\\
   \frac{\partial \J}{\partial t}&=-4 \eta \J +\tfrac{1}{2}
   (a-\lambda),
 \end{align*}
\end{rem}
To compute the curvature of $N^2$ we use the gauss equation
\eqref{gaussmetric} and obtain:
\begin{equation}\label{curvN2} 
   K(N^2)=e^{2f}(a-\J^2+2\M^2+\eta^2),
\end{equation}
which we verify by a straightforward computation is indeed independent
of $t$. 

Our first goal is to find out how $N^2$ is immersed in $\Rf$, i.~e. to
find an immersion independent of $t$. A look at the structure
equations \eqref{D11} - \eqref{D3xi} suggests to start with a linear
combination of $e_1$ and $\xi$.

We will solve the problem in two steps. First we define $v:=A e_1
+\xi$ for some function $A$ on $M^3$. Then $\pt v=\alpha v$ iff
$\alpha=A$ and $\pt A= A^2 -2\J A+ \lambda$, and $A:=-(\eta+\J)$
solves the latter differential equation. Next we define a positive
function $\beta$ on $\R$ as solution of the differential equation:
\begin{equation}\label{dtbeta}
\tfrac{\partial}{\partial t} \beta = \beta (\eta+\J)
\end{equation}
with initial condition $\beta(t_0)>0$. Then $\pt(\beta v)=0$ and by
\eqref{D21}, \eqref{D2xi}, \eqref{D31} and \eqref{D3xi} we get (since
$\beta$, $\eta$ and $\J$ only depend on $t$):
\begin{align}
D_{e_1}(\beta(-(\eta+\J) e_1 +\xi))&=0,\label{eq31}\\
D_{e_2}(\beta(-(\eta+\J) e_1 +\xi))&=-\beta(a +\eta^2 -\J^2)
e_2,\label{eq32}\\ 
D_{e_3}(\beta(-(\eta+\J)e_1 +\xi))&=-\beta(a +\eta^2 -\J^2)
e_3.\label{eq33} 
\end{align}

\begin{lem}\label{defnu}
  Define $\nu:=a +\eta^2 -\J^2$ on $\R$. $\nu$ is globally defined,
  $\pt(e^{2f} \nu)=0$ and $\nu$ vanishes identically or nowhere on $\R$.
\end{lem}
\begin{proof} Since $0=\pt K(N^2) = \pt(e^{2f}(\nu+2\M^2))$
  (cp.~\eqref{diffmu2} and \eqref{deff}) and $\pt(e^{2f}2 \M^2)=0$, we
  get that $\pt(e^{2f} \nu)=0$. Thus $\pt\nu=-2 (\pt f)\nu= -2\eta\nu$.
\end{proof}
Now we consider different cases depending on the behaviour of $\nu$.

\subsection{The first case: $\nu \neq 0$ on $M^3$}
\label{sec:case1}

We may, by translating $f$, i.e. by replacing $N^2$ with a homothetic
copy of itself, assume that $e^{2f} \nu =\epsilon$, where $\epsilon =\pm 1$.

\begin{lem}\label{defphi}
$\phi:=\beta (-(\eta +\J) e_1 +\xi)\colon N^2 \to \R^4$ defines a proper affine sphere in
a 3-dimensional linear subspace of $\R^4$. $\phi$ is part of a quadric
iff $\M =0$.
\end{lem}
\begin{proof}
By construction we ensured that $\pt \phi=0$ (cp. \eqref{eq31}), thus $\phi$
is defined on $N^2$. Furthermore it is an immersion, since $\phi_*(e_a)=
-\beta \nu e_a$ for $a= 2,3$ by \eqref{eq32} and \eqref{eq33}. A further
differentiation, using \eqref{D22} ($\beta$ and $\nu$ only depend on $t$),
gives:
\begin{align*}
D_{e_2} \phi_\star(e_2)& = -\beta (a +\eta^2-\J^2) D_{e_2}e_2\\
&= -\beta (a +\eta^2-\J^2) (-(\eta+\J) e_1 +\M e_2 + \ic223 e_3 +\xi)\\
&=\M\phi_\star(e_2)+\ic223 \phi_\star(e_3) -(a+\eta^2 -\J^2) \phi\\
&=\M\phi_\star(e_2)+\ic223 \phi_\star(e_3) -\epsilon e^{-2f} \phi.
\end{align*}
Similarly, we obtain the other derivatives,using \eqref{D23} -
\eqref{D33}, thus:
\begin{alignat*}{3}
D_{e_2} \phi_*(e_2)&= &\M \phi_*(e_2)  &+ \ic223 \phi_*(e_3)&- e^{-2f}\epsilon \phi
, \\
D_{e_2} \phi_*(e_3)&=& -\ic223 \phi_*(e_2)  &-\M \phi_*(e_3),& \\
D_{e_3} \phi_*(e_2)&=& &(\ic323-\M) \phi_*(e_3),& \\
D_{e_3} \phi_*(e_3)&=& -(\ic323+\M)\phi_*(e_2)  &&- e^{-2f}\epsilon \phi, \\
D_{e_a} \phi &=& -\beta e^{-2f}\epsilon e_a, & \qquad a=2,3.&
\end{alignat*}
We can read off the coefficients of the difference tensor $K^\phi$ of $\phi$
(cp. \eqref{strGauss} and \eqref{defK}): $(K^\phi)_{22}^2=\M$,
$(K^\phi)_{23}^3=-\M$, $(K^\phi)_{22}^3=0= (K^\phi)_{33}^3$, and see that
$\trace (K^\phi)_X$ vanishes. The affine metric introduced by this
immersion corresponds with the metric on $N^2$. Thus $-\epsilon \phi$ is the
affine normal of $\phi$ and $\phi$ is a proper affine sphere. Finally the
vanishing of the difference tensor characterizes quadrics.
\end{proof}

Our next goal is to find another linear combination of $e_1$ and $\xi$,
this time only depending on $t$. (Then we can express $e_1$ in terms
of $\phi$ and some function of $t$.)
\begin{lem}\label{defdelta}
  Define $\delta := a e_1 +(\eta -\J) \xi$. Then there exist a constant vector
  $C \in \R^4$ and a function $g(t)$ such that
$$ \delta(t)= g(t) C.$$
\end{lem}
\begin{proof} Using \eqref{D21} and \eqref{D2xi} resp. \eqref{D31} and
  \eqref{D3xi} we obtain that $D_{e_2}\delta = 0=D_{e_3} \delta$. Hence $\delta$
  depends only on the variable $t$. Moreover, we get by
  \eqref{e1a}, \eqref{D11}, \eqref{e1eta},\eqref{e1J} and \eqref{D1xi}
  that 
\begin{align*}
  \pt\delta&=D_{e_1} ( a e_1+(\eta-\J)\xi)\\
  &=(\J-\eta)(a-\lambda) e_1 + 2 a \J e_1+a \xi -(\eta-\J)\lambda e_1\\
  &\qquad +(-\eta^2 -3 \J^2 -a+ 4\J \eta)\xi\\
  &=(3 \J-\eta)(a e_1 +(\eta -\J)\xi)\\
  &=(3 \J -\eta) \delta.
\end{align*}
This implies that there exists a constant vector $C$ in $\mathbb
R^4$ and a function $g(t)$ such that $\delta(t)=g(t)C$.
\end{proof}
Combining $\phi$ and $\delta$ we obtain for $e_1$ (cp. Lem.~\ref{defphi} and
\ref{defdelta}) that
\begin{equation}\label{e1}
e_1(t,u,v)= -\tfrac{1}{\beta\nu}(\beta gC +(\eta-\J)\phi(u,v)).
\end{equation}
In the following we will use for the partial derivatives the
abbreviation $F_x:= \frac{\partial}{\partial x}F $, $x=t,u,v$.
\begin{lem}\label{partialF} 
\begin{align*}
&F_t = -\frac{g}{\nu}C -\pt(\frac{1}{\beta \nu})\phi,\\
&F_u = -\frac{1}{\beta\nu} \phi_u,\\
&F_v = -\frac{1}{\beta\nu} \phi_v.
\end{align*}
\end{lem}
\begin{proof} As by \eqref{dtbeta} and Lem.~\ref{defnu} $\pt
  \frac{1}{\beta\nu}= \frac{1}{\beta\nu}(\eta -\J)$, we obtain the equation for $F_t
  =e_1$ by \eqref{e1}. The other equations follow from \eqref{eq32}
  and \eqref{eq33}.
\end{proof}
It follows by the uniqueness theorem of first order differential
equations and applying a translation that we can write
$$F(t,u,v)= \tilde{g}(t) C -\frac{1}{\beta\nu}(t) \phi(u,v)$$
for a suitable
function $\tilde{g}$ depending only on the variable $t$. Since $C$ is
transversal to the image of $\phi$ (cp. Lem.~\ref{defphi} and
\ref{defdelta}), we obtain that after applying an equiaffine
transformation we can write: $F(t,u,v) =(\gamma_1(t), \gamma_2(t) \phi(u,v))$.
Thus we have proven the following:

\begin{thm} Let $M^3$ be an affine hypersurface of $\mathbb R^4$ which
  admits a pointwise $SO(2)$- or $\mathbb Z_3$-symmetry and with the
  globally defined function $(a +\eta^2 -\J^2)$ not identically zero
  on $M^3$. Then $M^3$ is affin equivalent to 
  $$F:I\times N^2\to \mathbb R^4:(t,u,v)\mapsto (\gamma_1(t),
  \gamma_2(t) \phi(u,v)),$$ 
  where
  $\phi: N^2 \to \mathbb R^3$ is an elliptic or hyperbolic affine sphere
  and $\gamma:I\to \mathbb R^2$ is an affine curve.\hfill \newline 
  Moreover,
  if $M^3$ admits a pointwise $SO(2)$-symmetry then $N^2$ is either an
  ellipsoid or a hyperboloid.
\end{thm}

In the next theorem we deal with the converse.
\begin{thm} Let $\phi:N^2 \to \mathbb R^3$ be an elliptic or hyperbolic
  affine sphere (scaled such that the absolute value of the mean
  curvature equals $1$) and let $\gamma: I \to \mathbb R^2$ be an affine
  curve.  Then if
  $$F(t,u,v)=(\gamma_1(t), \gamma_2(t) \phi(u,v)),$$
  defines a nondegenerate positive definite affine hypersurface, it
  admits a pointwise $\mathbb Z_3$- or $S_1$-symmetry.
\end{thm}

\begin{proof}
We have
\begin{align*}
  &F_t=(\gamma_1',\gamma_2' \phi),\\
  &F_u=(0,\gamma_2 \phi_u),\\
  &F_v=(0,\gamma_2 \phi_v),\\
  &F_{tt}=(\gamma_1'',\gamma_2'' \phi)=\tfrac{(\gamma_2''\gamma_1'-
  \gamma_1''\gamma_2')}{\gamma_1'}(0,\phi)+\tfrac{\gamma_1''}{\gamma_1'}
  F_t,\\  
  &F_{ut}=\tfrac{\gamma_2'}{\gamma_2} F_u,\\
  &F_{vt}=\tfrac{\gamma_2'}{\gamma_2} F_v,\\
  &F_{uu}=(0,\gamma_1\phi_{uu}),\\
  &F_{uv}=(0,\gamma_1\phi_{uv}),\\
  &F_{vv}=(0,\gamma_1\phi_{vv}).
\end{align*}
This implies that $F$ defines a nondegenerate affine immersion
provided that $\gamma_2\gamma_1 \gamma_1'(\gamma_2''\gamma_1'-\gamma_1''
\gamma_2') \neq 0$.  We moreover see that this immersion is definite
provided that the affine sphere is hyperbolic and $\gamma_1
\gamma_1'(\gamma_2''\gamma_1'-\gamma_1'' \gamma_2')>0$ or when the
proper affine sphere is elliptic and $\gamma_1
\gamma_1'(\gamma_2''\gamma_1'-\gamma_1'' \gamma_2')<0$.  As the proof
in both cases is similar, we will only treat the first case here. An
evaluation of the conditions for the affine normal $\xi$ ($\xi_t$,
$\xi_u$ , $\xi_v$ are tangential and $\det(F_t,F_u,F_v,\xi)=\sqrt{\det
  h}$) leads to:
\begin{equation*}
\xi=\alpha(t) (0,\phi(u,v)) +\beta(t) F_t,
\end{equation*}
where $(\gamma_2''\gamma_1'-\gamma_1'' \gamma_2')\gamma_1^2 =
\gamma_2^4 (\gamma_1')^3 \alpha^5$ and $\alpha'
+\tfrac{\beta(\gamma_2''\gamma_1'-\gamma_1''\gamma_2')}{\gamma_1'}=0$.
Taking $e_1$ in the direction of $F_t$, we see that $F_u$ and $F_v$
are orthogonal to $e_1$. It is also clear that $S$ restricted to the
space spanned by $F_u$ and $F_v$ is a multiple of the identity, and
$S(F_t)=\lambda F_t$, since $S$ is symmetric. Moreover, we have that
\begin{align*}
  (\nabla h)(F_t,F_u,F_u)
  &=(\tfrac{\gamma_1'}{\gamma_1}-\tfrac{\alpha'}{\alpha}-2
  \tfrac{\gamma_2'}{\gamma_2})h(F_u,F_u),\\
  (\nabla h)(F_t,F_u,F_v)
  &=(\tfrac{\gamma_1'}{\gamma_1}-\tfrac{\alpha'}{\alpha}-2
  \tfrac{\gamma_2'}{\gamma_2}) h(F_u,F_v),\\
  (\nabla h)(F_t,F_v,F_v)
  &=(\tfrac{\gamma_1'}{\gamma_1}-\tfrac{\alpha'}{\alpha}-2
  \tfrac{\gamma_2'}{\gamma_2}) h(F_v,F_v),\\
  (\nabla h)(F_u,F_t,F_t)&=0=(\nabla h)(F_v,F_t,F_t),
\end{align*}
implying that $K_{F_t}$ restricted to the space spanned by $F_u$ and
$F_v$ is a multiple of the identity.  Using the symmetries of $K$ it
now follows immediately that $F$ admits an $\mathbb Z_3$-symmetry or
an $SO(2)$- symmetry.
\end{proof}

\subsection{The second case: $\nu \equiv 0$ and $\J\neq \eta$ on $M^3$}
\label{sec:case2}

Next, we consider the case that $a =\J^2-\eta^2$ and $\eta \neq \J$ on
$M^3$. Since by \eqref{e1J} and \eqref{e1eta} $e_1(\eta -\J)=-\eta^2 -3 \J^2
-a+4\J \eta=4 \J (\eta -\J)$ we see that $\eta \neq \J$ everywhere on
$M^3$ or nowhere. 

We already have seen that $M^3$ admits a warped product structure. The
map $\phi$ we have constructed in Lemma~\ref{defphi} will not define
an immersion (cp. \eqref{eq32} and \eqref{eq33}). Anyhow, for a fixed
point $t_0$, we get from \eqref{D22} - \eqref{D33}, \eqref{eq32} and
\eqref{eq33}, using the notation $\tilde{\xi}=-(\eta+\J)
e_1 + \xi$:
\begin{align*}
&D_{e_2}e_2 =\ic223 e_3 +\M e_2 +\tilde{\xi},\\
&D_{e_2}e_3 =-\ic223 e_2 -\M e_3,\\
&D_{e_3}e_2 =\ic332 e_3 -\M e_3,\\
&D_{e_3}e_3 =\ic332 e_2 -\M e_2 +\tilde{\xi},\\
&D_{e_a}\tilde{\xi}=0, \qquad a=2,3.
\end{align*}

Thus, if $u$ and $v$ are local coordinates which span the second
distribution, then we can interprete $F(t_0,u,v)$ as an improper
affine sphere in a $3$-dimensional linear subspace.

Moreover, we see that this improper affine sphere is a paraboloid
provided that $\M$ at $t_0$ vanishes identically (as a function of $u$
and $v$). From the differential equations \eqref{diffmu2} determining
$\M$, we see that this is the case exactly when $\M$ vanishes
identically, i.e.  when $M$ admits a pointwise $SO(2)$-symmetry.

After applying a translation and a change of coordinates, we may
assume that
\begin{equation*}
F(t_0,u,v)=(u,v,f(u,v),0),
\end{equation*}
with affine normal $\tilde{\xi}(t_0,u,v)=(0,0,1,0)$. To obtain $e_1$
at $t_0$, we consider \eqref{D21} and \eqref{D31} and get that
\begin{equation*}
D_{e_a}(e_1 -(\eta -\J) F) = 0,\qquad a=2,3.
\end{equation*}
Evaluating at $t=t_0$, this means that there exists a constant vector
$C$ such that $e_1(t_0,u,v)=(\eta -\J)(t_0) F(t_0,u,v) +C$. Since
$\eta\neq J$ everywhere, we can write:
\begin{equation}\label{e1t0}
  e_1(t_0,u,v)=\alpha_1 (u,v, f(u,v),\alpha_2),
\end{equation}
where $\alpha_1\neq 0$ and we applied an equiaffine transformation so
that $C=(0,0,0,\alpha_1\alpha_2)$. To obtain information about $\pt
e_1$ we have that $D_{e_1}e_1= 2\J e_1 +\xi$ (cp. \eqref{D11}) and
$\xi=\tilde{\xi} + (\eta+\J)e_1$ by the definition of $\tilde{\xi}$.
Also we know that $\tilde{\xi}(t_0,u,v)=(0,0,1,0)$ and by \eqref{eq31}
- \eqref{eq33} that $D_{e_i}(\beta\tilde{\xi})=0$, $i=1,2,3$. Taking
suitable initial conditions for the function $\beta$ ($\beta(t_0)=1$),
we get that $\beta\tilde{\xi}=(0,0,1,0)$ and finally the following
vector valued differential equation:
\begin{equation}\label{dte1}
\pt e_1 =(\eta +3\J) e_1 +\beta^{-1} (0,0,1,0).
\end{equation}
Solving this differential equation, taking into account the initial
conditions \eqref{e1t0} at $t=t_0$, we get that there exist functions
$\delta_1$ and $\delta_2$ depending only on $t$ such that
\begin{equation*}
  e_1(t,u,v)= (\delta_1(t) u,\delta_1(t) v, \delta_1(t) (f(u,v)
  +\delta_2(t)), \alpha_2 \delta_1(t)), 
\end{equation*}
where $\delta_1(t_0)=\alpha_1$, $\delta_2(t_0)=0$, $\delta_1'(t)
=(\eta +3 \J) \delta_1(t)$ and $\delta_2'(t) =\delta_1^{-1}(t)
\beta^{-1}(t)$.  As $e_1(t,u,v) =\tfrac{\partial F}{\partial
  t}(t,u,v)$ and $F(t_0,u,v) =(u,v,f(u,v),0)$ it follows by
integration that
$$F(t,u,v)= (\gamma_1(t) u, \gamma_1(t) v, \gamma_1(t) f(u,v)
+\gamma_2(t) , \alpha_2 (\gamma_1(t)-1)),$$
where $\gamma_1'(t)
=\delta_1(t)$, $\gamma_1(t_0)=1$, $\gamma_2(t_0)=0$ and $\gamma_2'(t)
=\delta_1(t)\delta_2(t)$.  After applying an affine transformation
we have shown:

\begin{thm} \label{th3} Let $M^3$ be an affine hypersurface of
  $\mathbb R^4$ which admits a pointwise $SO(2)$- or $\mathbb
  Z_3$-symmetry and with the globally defined functions satisfying $a
  =-\eta^2 +\J^2$ but $a \not\equiv 0$ on $M^3$. Then $M^3$ is affine
  equivalent with
  $$F:I\times N^2\to \mathbb R^4:(t,u,v)\mapsto (\gamma_1(t) u,
  \gamma_1(t) v, \gamma_1(t) f(u,v) +\gamma_2(t),\gamma_1(t)),$$
  where
  $\psi: N^2 \to \mathbb R^3:(u,v) \mapsto (u,v,f(u,v))$ is an
  improper affine sphere with affine normal $(0,0,1)$ and $\gamma:I\to
  \mathbb R^2$ is an affine curve.\hfill \newline Moreover, if $M^3$
  admits a pointwise $SO(2)$-symmetry then $N^2$ is an elliptic
  paraboloid.
\end{thm}

In the next theorem we again deal with the converse.
\begin{thm} Let $\psi:N^2 \to \mathbb R^3:(u,v)\mapsto (u,v,f(u,v))$ be an
  improper affine sphere with affine normal
  $(0,0,1)$ and let $\gamma: I \to \mathbb R^2$ be an affine curve.
  Then if
  $$F(t,u,v)=(\gamma_1(t)u, \gamma_1(t)v,\gamma_1(t)f(u,v)
  +\gamma_2(t) ,\gamma_1(t)),$$ 
  defines a nondegenerate positive definite affine hypersurface, it
  admits a pointwise $\mathbb Z_3$- or $S_1$-symmetry.
\end{thm}
\begin{proof}
We have
\begin{align*}
  &F_t=(\gamma_1' u,\gamma_1' v ,\gamma_1' f(u,v) +\gamma_2',\gamma_1'),\\
  &F_u=(\gamma_1,0,\gamma_1 f_u,0),\\
  &F_v=(0,\gamma_1,\gamma_1 f_v,0),\\
  &F_{tt}=(\gamma_1'' u, \gamma_1''v ,\gamma_1'' f(u,v) +
  \gamma_2'',\gamma_1'')=\tfrac{\gamma_1''}{\gamma_1'}F_t
  +\tfrac{(\gamma_2''\gamma_1'-
    \gamma_1''\gamma_2')}{\gamma_1'}(0,0,1,0),\\
  &F_{ut}=\tfrac{\gamma_1'}{\gamma_1} F_u,\\
  &F_{vt}=\tfrac{\gamma_1'}{\gamma_1} F_v,\\
  &F_{uu}=(0,0,f_{uu} \gamma_1,0),\\
  &F_{uv}=(0,0,f_{uv} \gamma_1,0),\\
  &F_{vv}=(0,0,f_{vv}\gamma_1,0).
\end{align*}
This implies that $F$ defines a nondegenerate affine immersion
provided that $\gamma_1 \gamma_1'(\gamma_2''\gamma_1'-
\gamma_1''\gamma_2') \neq 0$.  We moreover see that this immersion is
definite provided that the improper affine sphere is positive definite
and $\gamma_1\gamma_1'(\gamma_2''\gamma_1'- \gamma_1''\gamma_2')>0$ or
when the improper affine sphere is negative definite and
$\gamma_1\gamma_1'(\gamma_2''\gamma_1'- \gamma_1''\gamma_2')<0$.  As
the proof in both cases is similar, we will only treat the first case
here. It easily follows that we can write the affine normal $\xi$ as:
\begin{equation*}
\xi=\alpha(t) (0,0,1,0) +\beta(t) F_t,
\end{equation*}
where $(\gamma_2''\gamma_1'-\gamma_1'' \gamma_2')=
\gamma_1^2(\gamma_1')^3 \alpha^5$ and $\alpha'
+\tfrac{\beta(\gamma_2''\gamma_1'-\gamma_1''\gamma_2')}{\gamma_1'}=0$.
Taking $e_1$ in the direction of $F_t$, we see that $F_u$ and $F_v$
are orthogonal to $e_1$. It is also clear that $S$ restricted to the
space spanned by $F_u$ and $F_v$ is a multiple of the identity, and
$S(F_t)=\lambda F_t$, since $S$ is symmetric. Moreover, we have that
\begin{align*}
  (\nabla h)(F_t,F_u,F_u)
  &=(-\tfrac{\gamma_1'}{\gamma_1}-\tfrac{\alpha'}{\alpha})
  h(F_u,F_u),\\
  (\nabla h)(F_t,F_u,F_v)
  &=(-\tfrac{\gamma_1'}{\gamma_1}-\tfrac{\alpha'}{\alpha})
  h(F_u,F_v),\\
  (\nabla h)(F_t,F_v,F_v)
  &=(-\tfrac{\gamma_1'}{\gamma_1}-\tfrac{\alpha'}{\alpha})
  h(F_v,F_v),\\
  (\nabla h)(F_u,F_t,F_t)&=0=(\nabla h)(F_v,F_t,F_t),
\end{align*}
implying that $K_{F_t}$ restricted to the space spanned by $F_u$ and
$F_v$ is a multiple of the identity.  Using the symmetries of $K$ it
now follows immediately that $F$ admits an $\mathbb Z_3$-symmetry or
an $SO(2)$- symmetry.
\end{proof}

\subsection{The third case: $\nu \equiv 0$ and $\J= \eta$ on $M^3$}
\label{sec:case3}

The final case now is that $a =\J^2-\eta^2$ and $\eta = \J$ on the
whole of $M^3$. This is dealt with in the following theorem:
\begin{thm} Let $M^3$ be an affine hypersurface of $\mathbb R^4$ which
  admits a pointwise $SO(2)$- or $\mathbb Z_3$-symmetry
  and with the globally defined functions satisfying $a =-\eta^2
  +\J^2$ and $\eta =\J$ on $M^3$. Then $M^3$ is affine equivalent to 
  $$F:I\times N^2\to \mathbb R^4:(t,u,v)\mapsto (u, v, f(u,v)
  +\gamma_2(t),\gamma_1(t)),$$
  where $\psi: N^2 \to \mathbb R^3:(u,v)
  \mapsto (u,v,f(u,v))$ is an improper affine sphere with affine
  normal $(0,0,1)$ and $\gamma:I\to \mathbb R^2$ is an affine
  curve.\hfill \newline Moreover, if $M^3$ admits a pointwise
  $SO(2)$-symmetry then $N^2$ is an elliptic paraboloid.
\end{thm}
\begin{proof} We proceed in the same way as in Theorem \ref{th3}. 
  We again use that $M^3$ admits a warped product structure and we fix a
  parameter $t_0$. At the point $t_0$, we have for $\tilde{\xi}=-(\eta+\J)
e_1 + \xi = -2\J e_1 + \xi$:
\begin{align*}
&D_{e_2}e_2 =\ic223 e_3 +\M e_2 +\tilde{\xi},\\
&D_{e_2}e_3 =-\ic223 e_2 -\M e_3,\\
&D_{e_3}e_2 =\ic332 e_3 -\M e_3,\\
&D_{e_3}e_3 =\ic332 e_2 -\M e_2 +\tilde{\xi},\\
&D_{e_a}\tilde{\xi}=0, \qquad a=2,3.
\end{align*}

Thus, if $u$ and $v$ are local coordinates which span the second
distribution, then we can interprete $F(t_0,u,v)$ as an improper
affine sphere in a $3$-dimensional linear subspace.

Moreover, we see that this improper affine sphere is a paraboloid
provided that $\M$ at $t_0$ vanishes identically (as a function of $u$
and $v$). From the differential equations \eqref{diffmu2} determining
$\M$, we see that this is the case exactly when $\M$ vanishes
identically, i.e.  when $M$ admits a pointwise $SO(2)$-symmetry.

After applying a translation and a change of coordinates, we may
assume that
\begin{equation*}
F(t_0,u,v)=(u,v,f(u,v),0),
\end{equation*}
with affine normal $\tilde{\xi}(t_0,u,v)=(0,0,1,0)$. To obtain $e_1$
at $t_0$, we consider \eqref{D21} and \eqref{D31} and get that
\begin{equation*}
D_{e_a}e_1 =(\eta -\J) e_a = 0,\qquad a=2,3.
\end{equation*}
It follows that $e_1(t_0,u,v)$ is a constant vector field. As it is
transversal, we may assume that there exists an $\alpha\neq 0$ such that
\begin{equation*}
e_1(t_0,u,v)= (0,0,0,\alpha).
\end{equation*}
As $e_1$ is determined by the differential equation
(cp. \eqref{dte1}):
\begin{equation*}
\frac{\partial e_1}{\partial t} =4\J  e_1 +\beta^{-1} (0,0,1,0),
\end{equation*}
it follows that 
\begin{equation*}
e_1(t,u,v)= (0,0,\delta_2(t),\delta_1(t)),
\end{equation*}
where $\delta_2(t_0)=0$, $\delta_2'(t)=4\J(t) \delta_2(t)+
\beta^{-1}(t)$, $\delta_1(t_0)=\alpha$ and $\delta_1'(t)=4\J(t)
\delta_1(t)$. Integrating once more with respect to $t$ we obtain
that
\begin{equation*}
F(t,u,v)= (u,v, f(u,v)+\gamma_2(t),\gamma_1(t)),
\end{equation*}
for some functions $\gamma_1$ and $\gamma_2$ with
$\gamma_i'=\delta_i$ and $\gamma_i(t_0)=0$ for $i=1,2$.
\end{proof}

In the next theorem we deal with the converse.

\begin{thm} Let $\psi:N^2 \to \mathbb R^3:(u,v)\mapsto (u,v,f(u,v))$ be an
  improper affine sphere with affine normal $(0,0,1)$ and let
  $\gamma: I \to \mathbb R^2$ be an affine curve.  Then if
  $$F(t,u,v)=(u, v,f(u,v) +\gamma_2(t) ,\gamma_1(t))$$
  defines a nondegenerate positive definite affine hypersurface, it
  admits a pointwise $\mathbb Z_3$- or $S_1$-symmetry.
\end{thm}

\begin{proof} We have
\begin{align*}
  &F_t=(0,0,\gamma_2',\gamma_1'),\\
  &F_u=(1,0,f_u,0),\\
  &F_v=(0,1,f_v,0),\\
  &F_{tt}=(0,0,\gamma_2'',\gamma_1'')=\tfrac{\gamma_1''}{\gamma_1'}F_t
  +\tfrac{(\gamma_2''\gamma_1'-
    \gamma_1''\gamma_2')}{\gamma_1'}(0,0,1,0),\\
  &F_{ut}=F_{vt}=0,\\
  &F_{uu}=(0,0,f_{uu},0),\\
  &F_{uv}=(0,0,f_{uv},0),\\
  &F_{vv}=(0,0,f_{vv},0).
\end{align*}
This implies that $F$ defines a nondegenerate affine immersion
provided that $\gamma_1'(\gamma_2''\gamma_1'-\gamma_1''\gamma_2') \neq
0$.  We moreover see that this immersion is definite provided that the
improper affine sphere is positive definite and
$(\gamma_2''\gamma_1'-\gamma_1''\gamma_2')\gamma_1'>0$ or when the
improper affine sphere is negative definite and
$(\gamma_2''\gamma_1'-\gamma_1''\gamma_2')\gamma_1'<0$.  As the proof
in both cases is similar, we will only treat the first case here. It
easily follows that we can write the affine normal $\xi$ as:
\begin{equation*}
\xi=\alpha(t) (0,0,1,0) +\beta(t) F_t,
\end{equation*}
where $(\gamma_2''\gamma_1'-\gamma_1''\gamma_2') = (\gamma_1')^3
\alpha^5$ and $\alpha'
+\tfrac{\beta(\gamma_2''\gamma_1'-\gamma_1''\gamma_2')}{\gamma_1'}=0$.
Taking now $e_1$ in the direction of $F_t$, we see that $F_u$ and
$F_v$ are orthogonal to $e_1$. It is also clear that $SF_u = SF_v=0$,
and $S(F_t)=\lambda F_t$, since $S$ is symmetric.  Moreover, we have
that
\begin{align*}
  (\nabla h)(F_t,F_u,F_u) &=-\tfrac{\alpha'}{\alpha} h(F_u,F_u)\\
  (\nabla h)(F_t,F_u,F_v) &=-\tfrac{\alpha'}{\alpha} h(F_u,F_v)\\
  (\nabla h)(F_t,F_v,F_v) &=-\tfrac{\alpha'}{\alpha} h(F_v,F_v),\\
  (\nabla h)(F_u,F_t,F_t)&=0=(\nabla h)(F_v,F_t,F_t),
\end{align*}
implying that $K_{F_t}$ restricted to the space spanned by $F_u$ and
$F_v$ is a multiple of the identity.  Using the symmetries of $K$ it
now follows immediately that $F$ admits an $\mathbb Z_3$-symmetry or
an $SO(2)$-symmetry.\end{proof}

\bigskip
\bigskip
\noindent
\begin{tabular}{l l}
Christine Scharlach  \hskip 10.5 cm &Luc Vrancken\\
Fakult{\"a}t II, Institut f\"ur Mathematik &Mathematisch Instituut\\
Technische Universit{\"a}t Berlin  & Universiteit Utrecht\\
Strasse des 17. Juni 136  & Budapestlaan 6\\
D-10623 Berlin& 3584CD Utrecht\\
Germany  & The Netherlands\\
cs\symbol{"40}math.tu-berlin.de & vrancken\symbol{"40}math.uu.nl
\end{tabular}


\begin{thebibliography}{10}
\bibitem{B}
R.~L.~Bryant,
\newblock {Second order families of special Lagrangian 3-folds,}
preprint, arXiv:math.DG/0007128.

\bibitem{D}
F.~Dillen, \newblock{ Equivalence theorems in affine differential geometry.}
\newblock{\em Geometriae Dedicata}, 
32 (1989):81--92.

\bibitem{DNV}
F.~Dillen, K.~Nomizu  and L.~Vrancken,
\newblock Conjugate connections and radon's theorem in affine differential
  geometry.
\newblock {\em Monatsh. Math.}, 109(1990): 221--235.




\bibitem{DV}
F.~Dillen and L.~Vrancken,
\newblock Calabi type composition of affine spheres.
\newblock {\em Differential Geometry and its applications}, 4(1994): 303--328.


\bibitem{H}
 S.~Hiepko,
 \newblock Eine innere Kennzeichung der verzerrten Produkte.
 \newblock {\em Math. Ann.}, 241(1979): 209--215.

 
\bibitem{LSZ} 
A.~M.~Li, U. Simon and G. Zhao, 
\newblock {\em Global
    Affine Differential Geometry of Hypersurfaces}, volume~11 of {\em
    De Gruyter Expositions in Mathematics}.  
\newblock Walter De Gruyter, Berlin-New York, 1993.

\bibitem{no}
S.~N{\"o}lker,
\newblock Isometric immersions of warped products.
\newblock{\em Differential Geometry and its applications}, 6(1996):1--30.

\bibitem{NS}
K.~Nomizu and T.~Sasaki,
\newblock {\em Affine Differential Geometry}.
\newblock Cambridge University press, Cambridge, 1994.


\bibitem{sv}
C.~Scharlach and L.~Vrancken,
\newblock Three dimensional affine hypersurfaces admitting pointwise
symmetries, 
preprint. 


\bibitem{vr00}
L.~Vrancken,
\newblock 
The Magid-Ryan conjecture for equiaffine hyperspheres with constant
  sectional curvature,
\newblock {\em Journal of Differential Geometry}, 54(2000): 99--138.

\bibitem{V}
L.~Vrancken,
\newblock{Special classes of three dimensional affine hyperpsheres
  characterized by properties of their cubic form,} preprint.


\end{thebibliography}
 \end{document}